\documentclass[11pt]{article}
\usepackage{amssymb,amsmath,color,graphicx} 

\def\R{\mathbb{R}}

\def\P{{\mathcal P}}
\def\A{{\bf A}}
\def\x{{\bf x}}
\def\b{{\bf b}}
\def\a{{\bf a}}
\def\1{{\bf 1}}
\def\bmatrix{\left( \begin{array}{ccccccccccccccccccccccccccc}} 
\def\ematrix{\end{array} \right)} 
\def\c{\cdots} 
\newcommand\degree{\operatorname{deg}} 

\definecolor{gray}{rgb}{.5,.5,.5} 
\def\gray{\color{gray}} 
\definecolor{black}{rgb}{1,1,1}

\begin{document}
\bibliographystyle{amsplain}

\newtheorem{theorem}{Theorem} 
\newtheorem{conjecture}{Conjecture} 

\begin{center} 

\Large{\bf The Number of ``Magic" Squares, Cubes, and Hypercubes} 
\footnote{Appeared in \emph{American Mathematical Monthly} {\bf 110}, no.~8 (2003), 707--717.} 
\normalsize

\vspace{12pt} 
{\sc Matthias Beck, Moshe Cohen, Jessica Cuomo, and Paul Gribelyuk} 

\end{center}



\section{INTRODUCTION.} 

\small
\begin{quote} 
The peculiar interest of magic squares and all \emph{lusus numerorum} in general lies in 
the fact that they possess the charm of mystery. They appear to betray some hidden 
intelligence which by a preconceived plan produces the impression of intentional design, 
a phenomenon which finds its close analogue in nature. \\ 
Paul Carus \cite[p.~vii]{andrewsws} 
\end{quote} 
\normalsize

Magic squares have turned up throughout history, some in a mathematical context, 
others in philosophical or religious contexts.  According to legend, the first magic square was 
discovered in China by an unknown mathematician sometime before the first century A.D. 
It was a magic square of order three thought to have 
appeared  on the back of a turtle emerging from a river.  Other magic 
squares surfaced at various places around the world in the centuries following their discovery. 
Some of the more interesting examples were recorded in Europe during the 1500s.  Cornelius 
Agrippa wrote \emph{De Occulta Philosophia} in 1510. In it he describes the 
spiritual powers of magic squares and produces some squares of orders from three up to nine. 
His work, although influential in the mathematical community, enjoyed only 
brief success, for the counter-reformation and the witch hunts of the Inquisition began soon thereafter: 
Agrippa himself was accused of being allied with the 
devil. Although this story seems outlandish now, we cannot 
ignore the strange mystical ties magic squares seem to have with the world 
and nature surrounding us, above and beyond their mathematical significance.  

Despite the fact that magic squares have been studied for a long time, they are still the 
subject of research projects. These include both mathematical-historical 
research, such as the discovery of unpublished magic squares of Benjamin Franklin \cite{paslesfranklin}, 
and pure mathematical research, much of which is connected with the algebraic and combinatorial 
geometry of polyhedra (see, for example, \cite{ahmeddeloerahemmecke}, \cite{beckpixton}, and \cite{mudgal}). 
Aside from mathematical research, 
magic squares naturally continue to be an excellent source of topics for 
``popular" mathematics books (see, for example, \cite{andrewsws} or \cite{pickoverzen}). 
In this paper we explore counting functions that are associated with magic squares. 

We define a \emph{semi-magic square} to be a square matrix whose entries are nonnegative integers and whose rows and columns 
(called \emph{lines} in this setting) sum 
to the same number. A \emph{magic square} is a semi-magic square whose main diagonals also add up to the line sum. 
A \emph{symmetric magic square} is a magic square that is a symmetric matrix. 
A \emph{pandiagonal magic square} is a semi-magic square whose 
diagonals parallel to the \emph{main} diagonal from the upper left to the lower 
right, wrapped around (i.e., continued to a duplicate square placed to the left or 
right of the given one), add up to the line sum. 
Figure \ref{examples} illustrates our various definitions. 
We caution the reader about clashing definitions in the literature. For example, 
some people would reserve the term ``magic square" for what we will call a 
\emph{traditional magic square}, a magic square of order $n$ whose entries are the integers $ 1 , 2 , \dots, n^2 $. 

\begin{figure}[htb] 
\begin{align*} 
\begin{array}{|c|c|c|} 
\hline 
 3 & 0 & 0 \\ 
\hline 
 0 & 1 & 2 \\ 
\hline 
 0 & 2 & 1 \\ 
\hline 
\end{array} 
\qquad 
\begin{array}{|c|c|c|} 
\hline 
 1 & 2 & 0 \\ 
\hline 
 0 & 1 & 2 \\ 
\hline 
 2 & 0 & 1 \\ 
\hline 
\end{array} 
\qquad 
\begin{array}{|c|c|c|c|} 
\hline 
 1 & 0 & 1 & 2 \\ 
\hline 
 0 & 1 & 2 & 1 \\ 
\hline 
 1 & 2 & 1 & 0 \\ 
\hline 
 2 & 1 & 0 & 1 \\ 
\hline 
\end{array} 
\end{align*} 
\caption{A semi-magic, a magic, and a symmetric pandiagonal magic square.} 
\label{examples} 
\end{figure} 

Our goal is to count these various types of squares. In the traditional case, 
this is in some sense not very interesting: for each order there is a fixed number 
of traditional magic squares. For example, there are 880 traditional 
$4 \times 4$ magic squares. The situation becomes more interesting if we drop 
the condition of traditionality and study the number of magic 
squares as a function of the line sum. 
We denote the total number of semi-magic, magic, symmetric magic, 
and pandiagonal magic squares of order $n$ and 
line sum $t$ by $ H_n (t), M_n (t), S_n (t)$, and $ P_n (t) $, 
respectively. 

We illustrate these notions for the case $ n=2 $, which is not very complicated: 
here a semi-magic square is determined once we know 
one entry (denoted by $*$ in Figure \ref{n=2}); a magic square has to have identical 
entries in each coefficient. 
\begin{figure}[htb] 
\begin{align*} 
\begin{array}{|c|c|} 
\hline 
 * & t-* \\ 
\hline 
 t-* & * \\ 
\hline 
\end{array} 
\qquad 
\begin{array}{|c|c|} 
\hline 
 t/2 & t/2 \\ 
\hline 
 t/2 & t/2 \\ 
\hline 
\end{array} 
\end{align*} 
\caption{Semi-magic and magic squares for $n=2$.} 
\label{n=2} 
\end{figure} 
Hence 
  \[ H_2 (t) = t+1 \ , \qquad M_2 (t) = S_2 (t) = P_2 (t) = \left\{ \begin{array}{ll} 1 & \mbox{ if $t$ is even, } \\ 
                                                                                      0 & \mbox{ if $t$ is odd. } \end{array} \right.  \] 
These easy results already hint at something: the counting function $ H_n 
$ is of a different character than the functions $ M_n $, $ S_n $, and $ P_n $. 

The oldest nontrivial results on this subject were first published in 1915: Macmahon \cite{macmahon} proved that 
  \[ H_3 (t) = 3 {{ t+3 }\choose{ 4 }} + {{ t+2 }\choose{ 2 }} \ , \qquad M_3 (t) = \left\{ \begin{array}{ll} \frac{2}{9} t^2 + \frac{2}{3} t + 1 & \mbox{ if } 3|t , \\ 
                                                                                                                                                0 & \mbox{ otherwise. } \end{array} \right. \] 
The first structural result for general $n$ was proved in 1973 by Ehrhart \cite{ehrhartmagic} and Stanley \cite{stanleymagic}. 
\begin{theorem}[Ehrhart, Stanley]\label{polyrec}. 
The function $ H_n (t) $ is a polynomial in $t$ of degree $ (n-1)^2 $ that satisfies the identities 
  \[ H_n (-n-t) = (-1)^{n-1} H_n (t) \ , \qquad \quad H_n (-1) = H_n (-2) = \dots = H_n (-n+1) = 0 \ . \] \end{theorem} 
The fact that $H_n$ is a polynomial of degree $(n-1)^2$ 
was conjectured earlier by Anand, Dumir, and Gupta \cite{ananddumirgupta}. 
An elementary proof of Theorem \ref{polyrec} can be found in \cite{spencermagic}. 
Stanley also proved analogous results for the counting functions for 
\emph{symmetric} semi-magic squares. 
In this paper, we establish analogues of these theorems for the other counting functions mentioned earlier. 

A \emph{quasi-polynomial} $Q$ of \emph{degree} $d$ is an expression of the form
  \[ Q(t) = c_{d}(t) \ t^{d} + \dots + c_{1}(t) \ t + c_{0}(t) \ , \]
where $c_0,c_1,\dots,c_d$ are periodic functions of $t$ and $c_d \not\equiv 0$. 
The least common multiple of the periods of the $c_j$ is called the \emph{period} of
$Q$. With this definition we can state the first of our two main theorems. 
\begin{theorem}\label{main} The functions $ M_n (t) $, $ S_n (t)$, and $ P_n (t) $ are
quasi-polynomials in $t$ of degrees 
$ n^2 - 2n - 1 $, $ n^2/2 - n/2 - 2 $, and $ n^2 - 3n + 2 $, respectively, that satisfy the identities 
  \begin{equation}\label{reclaws} \begin{array}{ll} M_n (-n-t) = (-1)^{ n - 1 } \ M_n (t) \ , & \qquad M_n (-1) = M_n (-2) = \dots = M_n (-n+1) = 0 \ , \\ 
                                                    S_n (-n-t) = (-1)^{ n(n-1)/2 } \ S_n (t) \ , & \qquad S_n (-1) = S_n (-2) = \dots = S_n (-n+1) = 0 \ , \\ 
                                                    P_n (-n-t) = P_n (t) \ ,                  & \qquad P_n (-1) = P_n (-2) = \dots = P_n (-n+1) = 0 \ . \end{array} \end{equation}
\end{theorem} 

We can extend the counting function $H_n$ for semi-magic squares in 
the following way. Define a \emph{semi-magic hypercube} (also called a 
\emph{quasi-magic hypercube}) to be a $d$-dimensional 
$n \times \dots \times n$ array of $n^d$ nonnegative integers 
that sum to the same number $t$ parallel to any axis; that is, if we 
denote the entries of the array by $m_{j_1 \dots j_d} $ ($ 1 \leq j_k \leq n$), 
then we require $\sum_{j_k=1}^n m_{j_1 \dots j_d} = t$ for all $k=1, \dots, d$. 
Again we count all such cubes in terms of $d$, $n$, and $t$; 
we denote the corresponding enumerating function by $H_n^d (t)$. So 
$H_n^2=H_n$, whose properties are stated in Theorem \ref{polyrec}. 
Except for the case $d=n=3$, which seems to have appeared first in 
\cite{bonamagiccubes}, we could not find any other references. We prove the 
following result for general $d$. 
\begin{theorem}\label{cubes} The function $ H_n^d (t) $ is a quasi-polynomial in $t$ of 
degree $(n-1)^d$ that satisfies the identities 
  \[ H_n^d (-n-t) = (-1)^{ n - 1 } \ H_n^d (t) \ , \qquad \qquad H_n^d (-1) = H_n^d (-2) = \dots = H_n^d (-n+1) = 0 \ . \] 
\end{theorem}

We prove Theorem \ref{main} in Sections \ref{comb} and \ref{degrees}.
To be able to compute the counting functions $ M_n $, $ S_n $, and $ P_n $ for specific $n$, 
we need the periods of these quasi-polynomials. We describe in Section 
\ref{comp} methods for finding these periods and hence for the actual 
computation of $ M_n $, $ S_n $, and $ P_n $. Finally, we prove Theorem \ref{cubes} in Section 
\ref{hypercubes}. All of our methods are based on the idea that one can interpret the 
various counting functions as enumerating integer points (``lattice 
points'') in certain polytopes. 


\section{SOME GEOMETRIC COMBINATORICS.}\label{comb} 

A \emph{convex polytope} $\P$ in $\R^d$ is the convex hull of finitely many points 
in $\R^d$. Alternatively (and this correspondence is nontrivial \cite{ziegler}), one can 
define $\P$ as the bounded intersection of affine halfspaces. 
If there is a hyperplane $H = \left\{ \x \in \R^d : \ \a \cdot \x = b \right\} $ such that $\P \cap H$ 
consists of a single point, then this point is a \emph{vertex} of $\P$. 
A polytope is \emph{rational} if all of its vertices have rational coordinates. 

\begin{figure}[ht]
\begin{center}
\includegraphics{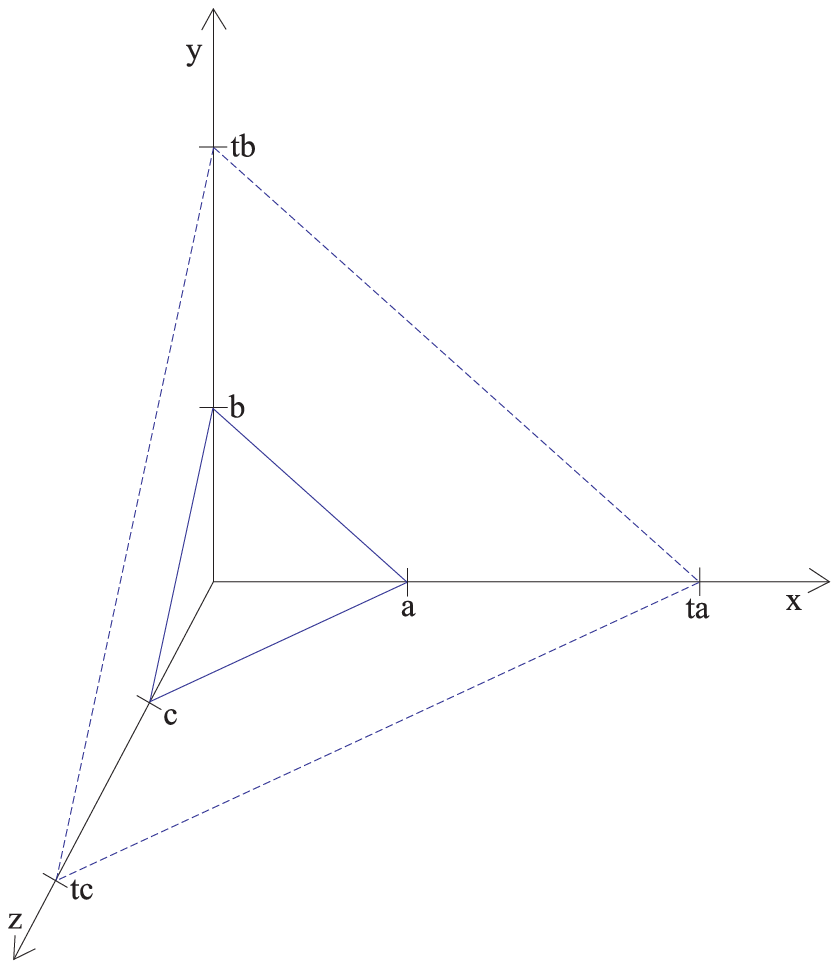} 
\end{center}
\caption{Dilation of a polytope.}\label{dilationfig}
\end{figure}

Suppose that $\P$ is a convex rational polytope in $\R^d$. For a positive integer $t$, let $ L_\P (t) $ denote the number 
of lattice points in the dilated polytope $ t \P = \{ tx : x \in \P  \} $ (see Figure \ref{dilationfig}). Similarly, define $ L_\P^{*} (t) $ to be the number of 
lattice points in the (relative) interior of $ t \P $. Ehrhart, who initiated the study of the lattice point count in dilated
polytopes, proved the following \cite{ehrhart2}. 
  \begin{theorem}[Ehrhart]\label{quasi}. If $\P$ is a convex rational polytope, then the functions 
  $ L_\P (t) $ and $ L_\P^{*} (t) $ are quasi-polynomials in $t$ whose degree is the dimension of $\P$ 
  and whose period divides the least common multiple of the denominators of the vertices of $\P$. \end{theorem} 
In particular, if $\P$ has integer vertices, then $ L_\P $ and $ L_\P^{*} $ are polynomials. 
Ehrhart conjectured and partially proved the following reciprocity law, which was proved by 
Macdonald \cite{macdonald} (for the case that $\P$ has integer vertices) 
and McMullen \cite{mcmullenreciprocity} (for the case that $\P$ has arbitrary rational vertices). 
   \begin{theorem}[Ehrhart-Macdonald reciprocity law]\label{rec}. 
     $ L_\P (-t) = (-1)^{\dim \P} L_\P^{*} (t) $ 
   \end{theorem} 

Theorems \ref{quasi} and \ref{rec} mark the beginning of our journey towards a proof of Theorem \ref{main}. 
The connection to our counting functions $ M_n $, $ S_n $, and $ P_n $ is the following: all the various magic square conditions are linear inequalities (the entries 
are nonnegative) and linear equalities (the entries in each line sum up to $t$) in the $ n^2 $ variables forming the entries of the 
square. In other words, what we are computing is the number of nonnegative 
integer solutions to the linear system 
  \[ \A \ \x = \b \ , \] 
where $ \x $ is a column vector in $ \R^{ n^2 } $ whose components are the entries of 
our square, $ \A $ denotes the matrix determining the 
magic-sum conditions, and $ \b $ is the column vector each of whose entries is $t$. 
We will use the convention that the entries of $\x$ are ordered by rows of the original 
square, and from left to right in each row, that is, the entries $x_{jk}$ of the square 
are being arranged as 
  \[ \x = \left( x_{11}, x_{12}, \dots, x_{1n}, x_{21}, x_{22}, \dots, x_{2n}, \dots, x_{n1}, x_{n2}, \dots, x_{nn} \right) \ . \] 
The matrix $\A$ is constructed accordingly. 
For example, if we study magic $3 \times 3$ squares, 
  \[ \A = \left( \begin{array}{ccccccccc} 1 & 1 & 1 & 0 & 0 & 0 & 0 & 0 & 0 \\ 
                                          0 & 0 & 0 & 1 & 1 & 1 & 0 & 0 & 0 \\ 
					  0 & 0 & 0 & 0 & 0 & 0 & 1 & 1 & 1 \\ 
					  1 & 0 & 0 & 1 & 0 & 0 & 1 & 0 & 0 \\ 
					  0 & 1 & 0 & 0 & 1 & 0 & 0 & 1 & 0 \\ 
					  0 & 0 & 1 & 0 & 0 & 1 & 0 & 0 & 1 \\ 
					  1 & 0 & 0 & 0 & 1 & 0 & 0 & 0 & 1 \\ 
					  0 & 0 & 1 & 0 & 1 & 0 & 1 & 0 & 0 \end{array} \right) \ . \] 
The first row of $\A \x = \b$ represents the first row sum in the square: 
$x_{11} + x_{12} + x_{13} = t$; the second and third row ($x_{21} + x_{22} + x_{23} = t$ 
and $x_{31} + x_{32} + x_{33} = t$) are the sums of the second and third row in the square. 
The fourth row of $\A \x = \b$ represents the first 
column sum in the square: $x_{11} + x_{21} + x_{31} = t$; and the next two rows take care of the 
second and third column in our square. Finally the last two rows of $\A \x = \b$ 
represent the two diagonal constraints $x_{11} + x_{22} + x_{33} = t$ and 
$x_{13} + x_{22} + x_{31} = t$. 

Furthermore, by writing $ \b = t \ \1 $, where $ \1 $ signifies a column vector each of whose entries is 1, we can see that our counting function 
enumerates nonnegative integer solutions to the linear system $ \A \ \x = t \ \1 $, that is, it counts the lattice points in the $t$-dilate of the polytope 
  \[ \P = \left\{ \x = \left( x_1 , \dots , x_{n^2} \right) \in \R^{n^2} : \ x_k \geq 0 \text{ for } k = 1, 2, \dots, n^2 , \ \A \ \x = \1 \right\} \ , \] 
provided that we choose the matrix $\A$ according to the magic-sum conditions. 
Note that $\P$ is an intersection of half-spaces and hyperplanes, and is therefore convex. 
No matter which counting function the matrix $\A$ corresponds to, the entries of $\A$ are 0s and 1s. We obtain the 
vertices of $\P$ by converting some of the inequalities $ x_k \geq 0 $ to equalities. It is easy to conclude from this that the vertices 
of $\P$ are rational. Hence by Theorem \ref{quasi}, 
$ M_n$, $S_n$, and $ P_n $ are quasi-polynomials whose degrees are the dimensions of the corresponding polytopes. 
Geometrically, the ``magic-sum variable'' $t$ is the dilation factor of these polytopes. 

The Ehrhart-Macdonald reciprocity law (Theorem \ref{rec}) connects the lattice-point count in $ \P $ to that of the interior of $\P$. In our 
case, this interior (again, using any matrix $\A$ suitable for one of the counting functions) is described by 
  \[ \P^{ \text{int} } = \left\{ \x = \left( x_1 , \dots , x_{n^2} \right) \in \R^{n^2} : \ x_k > 0 \text{ for } k = 1, 2, \dots, n^2 , \ \A \ \x = \1 \right\} \ . \] 
The lattice-point count is the same as before, with the difference that we now allow only {\it positive} integers as solutions to the linear system. 
This motivates us to define 
by $ M_n^{*} (t), S_n^{*} (t)$, and $ P_n^{*} (t) $ the counting functions for magic squares, symmetric magic squares, and pandiagonal 
magic squares as before, but with the restriction that the entries be positive integers. 
Theorem \ref{rec} then asserts, for example, that 
  \begin{equation}\label{law} M_n (-t) = (-1)^{ \degree (M_n) } M_n^{*} (t) \ . \end{equation} 
On the other hand, by its very definition $ M_n^{*} (t) = 0 $ for $ t = 1, 2, \dots, n-1 $. Hence we obtain $ M_n (-t) = 0 $ for such $t$. 
Also, since each row of the matrix $ \A $ defined by some magic-sum conditions has exactly $n$ 1s and all other entries 0, it is not hard to conclude 
that $ M_n^{*} (t) = M_n (t-n) $. Combining this with the reciprocity law (\ref{law}), we obtain 
  \[ M_n (t) = (-1)^{ \degree (M_n) } M_n (-t-n) \ . \] 
There are analogous statements for $ S_n $ and $ P_n $. 
Once we verify the degree formulas of Theorem \ref{main}, 
(\ref{reclaws}) follows from the observation that 
  \begin{eqnarray*} &\mbox{}& (-1)^{ n^2 - 2n - 1 } = (-1)^{ n - 1 } \ , \\ 
                    &\mbox{}& (-1)^{ \frac{1}{2} n^2 - \frac{1}{2} n - 2 } = (-1)^{ n(n-1)/2 } \ , \\ 
                    &\mbox{}& (-1)^{ n^2 - 3n + 2 } = (-1)^{ n(n-3) } = 1 \ . \end{eqnarray*} 


\section{PROOF OF THE DEGREE FORMULAS.}\label{degrees} 
Let's start with magic squares: by Theorem \ref{quasi}, the degree of $M_n$ is $n^2$ (the number of places we have to fill) minus 
the number of linearly independent constraints. In other words, we have to find the rank of the $n^2 \times (2n+2)$ matrix 
  \[ \A = \bmatrix 1 &    &    & \gray \mid  & 1 &    & & \gray \mid  &  & \gray \mid  & 1 \\[-5.5pt] 
	        & \ddots & & \gray \mid  & & \ddots & & \gray \mid  & \c & \gray \mid  &  & \ddots \\ 
		&    &  1 & \gray \mid  &   &    & 1  & \gray \mid  &  & \gray \mid  & &    & 1 \\[-4pt] 
              \gray - & \gray - & \gray - & \gray - & \gray - & \gray - & \gray - & \gray - & \gray - & \gray - & \gray - & \gray - & \gray - \\[-4pt] 
              1 & \c & 1  & \gray \mid  &   &    &  & \gray \mid  &   & \gray \mid  &  &   &    &   \\ 
                &    &    & \gray \mid  & 1 & \c & 1  & \gray \mid  &  &  \gray \mid  &  &    &    &   \\[-5.5pt] 
                &    &    & \gray \mid  &   &    &    & \gray \mid  & \ddots & \gray \mid  \\ 
                &    &    & \gray \mid  &   & & &  \gray \mid  &    &  \gray \mid  & 1 & \c & 1 \\[-4pt] 
              \gray - & \gray - & \gray - & \gray - & \gray - & \gray - & \gray - & \gray - & \gray - & \gray - & \gray - & \gray - & \gray - \\[-4pt] 
              1 &    &    & \gray \mid  &   & \hspace{-.3in} 1 &    & \gray \mid  & \c &  \gray \mid  &  &   & 1 \\ 
                &    &  1 & \gray \mid  &   &   1 \hspace{-.3in} &  &  \gray \mid  & \c & \gray \mid  & 1 & 
	     \ematrix \ . \] 
Here we only show the entries that are 1; every other entry is 0. 
Similarly  as in our example, the first $n$ rows of $\A$ represent the 
$n$ column constraints of our square ($x_{11} + \dots + x_{n1} = t, \ \dots, \ x_{1n} + \dots + x_{nn} = t$), 
the next $n$ rows represent the row contraints of the square, and the last two 
rows take care of the diagonal constraints. 

The sum of the first $n$ rows of $\A$ is the same as the sum of the next $n$ rows, so 
one of the first $2n$ rows is redundant and we can eliminate it; we choose the $(n+1)$th. Furthermore, we can add 
the difference 
of the first and $(n+2)$th row to the $(2n+1)$th and subtract the $n$th row
from the $(2n+2)$th. These operations yield the $n^2 \times (2n+1)$ matrix 
  \[ \bmatrix 1 &    &   & \gray \mid  & 1 &    &    &  \gray \mid  &   & \gray \mid  & 1 \\[-5.5pt] 
	        & \ddots & & \gray \mid & & \ddots & & \gray \mid  & \c & \gray \mid  &  & \ddots \\ 
		&    &  1  & \gray \mid  & &    & 1  & \gray \mid  &    & \gray \mid  &   &    & 1 \\[-4pt] 
              \gray - & \gray - & \gray - & \gray - & \gray - & \gray - & \gray - & \gray - & \gray - & \gray - & \gray - & \gray - & \gray - \\[-4pt] 
                &    &  & \gray \mid & 1 & \c & 1  & \gray \mid  &    & \gray \mid  &   &    &    &   \\[-5.5pt] 
                &    &    & \gray \mid  &    &    & & \gray \mid  & \ddots & \gray \mid  & \\ 
                &    &    & \gray \mid  &    &    &    & \gray \mid  & & \gray \mid  & 1 & \c & 1 \\[-4pt] 
              \gray - & \gray - & \gray - & \gray - & \gray - & \gray - & \gray - & \gray - & \gray - & \gray - & \gray - & \gray - & \gray - \\[-4pt] 
                &    &    & \gray \mid  & 0 & \hspace{-.09in} 2 \ 1 \ \c & 1 & \gray \mid  & * & \gray \mid & & *  \\ 
		&    &    & \gray \mid  &    & 1 \hspace{-.5in} & -1 & \gray \mid  & * & \gray \mid &  &  * 
  	     \ematrix \ . \] 
These represent $2n+1$ linearly independent restrictions on our magic square. 
The polytope corresponding to $M_n$ therefore has dimension $ n^2 - 2n - 1 $: it lies in an affine subspace 
of that dimension, and the point $ \left( 1/n , \dots , 1/n \right) $ is an interior 
point of the polytope. Hence, by Theorem \ref{quasi}, the degree of $M_n$ is $ n^2 - 2n - 1 $. 

The case of symmetric magic squares is simpler: here we have $ \sum_{k=1}^n k = n(n+1)/2 $ places 
to fill, and the row-sum condition is equivalent to the column-sum condition. Hence there are $ n+2 $ constraints, 
which are easily identified as linearly independent. The dimension of the 
corresponding polytope, and the degree of $S_n$, 
is therefore $ n^2/2 - n/2 - 2 $. 

Finally, we discuss pandiagonal magic squares. Again we have $n^2$ places to fill; this time the constraints are 
represented by the $n^2 \times 3n$ matrix 
  \[ \bmatrix 1 &    &    &   & \gray \mid  & 1  &    &    & & \gray \mid  &  & \gray \mid  & 1  \\ 
	        &    & \hspace{-.2in} \ddots & & \gray \mid  &  &    & \hspace{-.2in} \ddots & &  \gray \mid  & \c & \gray \mid  &  &    & \hspace{-.2in} \ddots \\ 
		&    &    & 1 & \gray \mid  &    &    &    & 1 & \gray \mid  &    & \gray \mid  &    &    &    & 1 \\[-4pt] 
              \gray - & \gray - & \gray - & \gray - & \gray - & \gray - & \gray - & \gray - & \gray - & \gray - & \gray - & \gray - & \gray - & \gray - & \gray - & \gray - \\[-4pt] 
              1 &    &    &   & \gray \mid  &    & 1  &    &   & \gray \mid  &    & \gray \mid  &    &    &    & 1 \\[-5.5pt] 
                & 1  &    &   & \gray \mid  &    &    & \ddots & & \gray \mid  & \c & \gray \mid  & 1 \\[-5.5pt] 
	        &    & \ddots & & \gray \mid  &  &    &    & 1 & \gray \mid  &    &  \gray \mid  &   & \ddots \\ 
                &    &    & 1 & \gray \mid  & 1  &    &    &   & \gray \mid  &    & \gray \mid  &    &    & 1  \\[-4pt] 
              \gray - & \gray - & \gray - & \gray - & \gray - & \gray - & \gray - & \gray - & \gray - & \gray - & \gray - & \gray - & \gray - & \gray - & \gray - & \gray - \\[-4pt] 
              1 &  & \hspace{-.2in} \c & 1 & \gray \mid  &    &    &    &   & \gray \mid  & & \gray \mid  & \\ 
                &    &    &   & \gray \mid  & 1  &  & \hspace{-.2in} \c & 1 & \gray \mid  & & \gray \mid  & \\[-5.5pt] 
                &    &    &   & \gray \mid  &    &    &    &   & \gray \mid  & \ddots & \gray \mid  & \\ 
                &    &    &   & \gray \mid  &    &    &    &   & \gray \mid  &    & \gray \mid  &  1 &  & \hspace{-.2in} \c & 1 \\ 
	     \ematrix \ . \] 
Here the first $n$ rows represent the column constraints in our square, the next 
$n$ rows the pandiagonal constraints, and the last $n$ rows the row constraints. 

As in the first case, the sum of the first $n$ rows equals the sum of the 
next $n$ rows, as well as the sum of the last $n$ rows. We can therefore eliminate 
the $(2n)$th and $(2n+1)$th rows to get the $n^2 \times (3n-2)$ matrix 
  \[ \bmatrix 1 &    &    &   & \gray \mid  & 1  &    &    & & \gray \mid  &  & \gray \mid  & 1  \\ 
	        &    & \hspace{-.2in} \ddots & & \gray \mid  &  &    & \hspace{-.2in} \ddots & &  \gray \mid  & \c & \gray \mid  &  &    & \hspace{-.2in} \ddots \\ 
		&    &    & 1 & \gray \mid  &    &    &    & 1 & \gray \mid  &    & \gray \mid  &    &    &    & 1 \\[-4pt] 
              \gray - & \gray - & \gray - & \gray - & \gray - & \gray - & \gray - & \gray - & \gray - & \gray - & \gray - & \gray - & \gray - & \gray - & \gray - & \gray - \\[-4pt] 
              1 &    &    &   & \gray \mid  &    & 1  &    &   & \gray \mid  &    & \gray \mid  &    &    &    & 1 \\[-5.5pt] 
                & \ddots &    &   & \gray \mid  &    &    & \ddots & & \gray \mid  & \c & \gray \mid  & 1 \\[-5.5pt] 
	        &    & 1 & & \gray \mid  &  &    &    & 1 & \gray \mid  &    &  \gray \mid  &   & \ddots & 1 \\[-4pt] 
              \gray - & \gray - & \gray - & \gray - & \gray - & \gray - & \gray - & \gray - & \gray - & \gray - & \gray - & \gray - & \gray - & \gray - & \gray - & \gray - \\[-4pt] 
                &    &    &   & \gray \mid  & 1  &  & \hspace{-.2in} \c & 1 & \gray \mid  & & \gray \mid  & \\[-5.5pt] 
                &    &    &   & \gray \mid  &    &    &    &   & \gray \mid  & \ddots & \gray \mid  & \\ 
                &    &    &   & \gray \mid  &    &    &    &   & \gray \mid  &    & \gray \mid  &  1 &  & \hspace{-.2in} \c & 1 \\ 
	     \ematrix \ . \] 
In this new matrix, we can now replace the $(n+1)$th row by the difference of the $(n+1)$th and first rows, 
the $(n+2)$th row by the difference of the $(n+2)$th and second rows, and so on:
  \[ \bmatrix 1 &     &   & \gray \mid  & 1  &        & & \gray \mid  &  & \gray \mid  & 1  \\[-5.5pt] 
	        & \ddots & & \gray \mid &    & \ddots & &  \gray \mid  & \c & \gray \mid  &  &   \ddots \\ 
		&     & 1 & \gray \mid  &    &      & 1 & \gray \mid  &    & \gray \mid  &    &     & 1 \\[-4pt] 
              \gray - & \gray - & \gray - & \gray - & \gray - & \gray - & \gray - & \gray - & \gray - & \gray - & \gray - & \gray - & \gray - \\[-4pt] 
                &     &   & \gray \mid  & -1 & \hspace{-.3in} 1 & & \gray \mid  &    & \gray \mid  &  \\[-5.5pt] 
                &     &   & \gray \mid  &    & \ddots           & & \gray \mid  & * & \gray \mid  & & * \\ 
	        &     &   & \gray \mid  &  &   -1 \hspace{-.2in}  & 1 & \gray \mid  &    &  \gray \mid  &    \\[-4pt] 
              \gray - & \gray - & \gray - & \gray - & \gray - & \gray - & \gray - & \gray - & \gray - & \gray - & \gray - & \gray - & \gray - \\[-4pt] 
                &    &    &   \gray \mid  & 1  &  \c & 1 & \gray \mid  & & \gray \mid  & \\[-5.5pt] 
                &    &    &   \gray \mid  &    &    &    &  \gray \mid  & \ddots & \gray \mid  & \\ 
                &    &    &   \gray \mid  &    &    &    &  \gray \mid  &    & \gray \mid  &  1 &  \c & 1 \\ 
	     \ematrix \ . \] 
Finally, we can replace the $(n+2)$th row by the sum of the $(n+1)$th and $(n+2)$th rows,
the $(n+3)$th row by the sum of the $(n+2)$th and $(n+3)$th rows, etc. Subtracting
the $(2n)$th row from the $(n+1)$th row gives a matrix of full rank, that is, rank $3n-2$.
The dimension of the corresponding polytope, which is the degree of $P_n$, is thus $ n^2 - 3n + 2 $. 

This finishes the proof of Theorem \ref{main}. 


\section{COMPUTATIONS.}\label{comp} 
To interpolate a quasi-polynomial of degree $d$ and period $p$, we need to compute $p (d+1)$ values. 
The periods of $M_n$, $S_n$, and $P_n$ are not as simple to derive as their degrees. What we can do, however, is to 
compute for fixed $n$ the vertices of the respective polytopes, whose denominators give the periods of the quasi-polynomials (Theorem \ref{quasi}). 
This is easy for very small $n$ but gets complicated very quickly. For example, we can 
practically find by hand that the vertices of the polytope corresponding to $S_3$ are 
$\left( 2/3 , 0 , 1/3 , 1/3 , 2/3 , 0 \right) $ and $\left( 0 , 2/3 , 1/3 , 1/3 , 0 , 2/3  \right) $, 
whereas the polytope corresponding to $S_5$ has seventy-four vertices. To make computational matters worse, the 
least common multiple of the denominators of these vertices is 60; one would have to do a lot of interpolation to obtain $S_5$. 

The reciprocity laws in Theorem \ref{main} essentially halve the number of computations 
for each quasi-polynomial. Nevertheless, the task of computing our quasi-polynomials becomes impractical without 
further tricks or large computing power. The following table contains data about the vertices of the polytopes 
corresponding to $M_n$, $S_n$, and $P_n$ for $n=3$ and $4$. It was produced using {\tt Maple}. 

\begin{center} 
\begin{tabular}{c|c|c} Polytope corresponding to & Number of vertices & l.c.m.~of denominators \\ \hline 
                      $M_3$ &      4        & 3 \\ 
                      $S_3$ &      2        & 3 \\ 
                      $P_3$ &      3        & 1 \\ 
                      $M_4$ &     20        & 2 \\ 
                      $S_4$ &     12        & 4 \\ 
                      $P_4$ &     28        & 2 \end{tabular} \end{center} 

With this information, it is now easy (that is, easy with the aid of a computer) to interpolate each quasi-polynomial. 
Here are the results: 
\begin{eqnarray*} M_3 (t) &=& \left\{ \begin{array}{ll} \tfrac{2}{9} t^2 + \tfrac{2}{3} t + 1 & \mbox{ if } 3|t , \\[5pt]
                                                                         \quad 0 & \mbox{ otherwise; } \end{array} \right. \\[10pt] 
     S_3 (t) &=& \left\{ \begin{array}{ll} \tfrac{2}{3} t + 1 & \mbox{ if } 3|t , \\[5pt]
                                                      \quad 0 & \mbox{ otherwise; } \end{array} \right. \\[10pt]
     P_3 (t) &=& \tfrac{1}{2} t^2 + \tfrac{2}{3} + 1 ; \\[10pt]
     M_4 (t) &=& \left\{ \begin{array}{ll} \tfrac{1}{480}t^7 + \tfrac{7}{240}t^6 + \tfrac{89}{480}t^5 + \tfrac{11}{16}t^4 + \tfrac{49}{30}t^3 + \tfrac{38}{15}t^2 + \tfrac{71}{30}t + 1 & \mbox{ if } t \mbox{ is even, } \\[5pt] 
                                           \tfrac{1}{480}t^7 + \tfrac{7}{240}t^6 + \tfrac{89}{480}t^5 + \tfrac{11}{16}t^4 + \tfrac{779}{480}t^3 + \tfrac{593}{240}t^2 + \tfrac{1051}{480}t + \tfrac{13}{16} & \mbox{ if } t \mbox{ is odd; } \end{array} \right. \\[10pt] 
     S_4 (t) &=& \left\{ \begin{array}{ll} \tfrac{5}{128} t^4 + \tfrac{5}{16} t^3 + t^2 + \tfrac{3}{2} t + 1 & \mbox{ if } t \equiv 0 \mod 4 , \\[5pt] 
                                           \tfrac{5}{128} t^4 + \tfrac{5}{16} t^3 + t^2 + \tfrac{3}{2} t + \tfrac{7}{8} & \mbox{ if } t \equiv 2 \mod 4 , \\[5pt] 
                                           0 & \mbox{ if } t \mbox{ is odd; } \end{array} \right. \\[10pt] 
     P_4 (t) &=& \left\{ \begin{array}{ll} \tfrac{7}{1440}t^6 + \tfrac{7}{120}t^5 + \tfrac{23}{72}t^4 + t^3 + \tfrac{341}{180}t^2 + \tfrac{31}{15}t + 1 & \mbox{ if } t \mbox{ is even, } \\[5pt] 
                                           \quad 0 & \mbox{ if } t \mbox{ is odd. } \end{array} \right. \end{eqnarray*} 


\section{SEMI-MAGIC HYPERCUBES.}\label{hypercubes} 
We now prove Theorem \ref{cubes}. 
The fact that $ H_n^d $ is a quasi-polynomial, the reciprocity law for this function and $ H_n^{d*} $, 
and the location of special zeros for $ H_n^d $ 
follow in exactly the same way as the respective statements in Theorem \ref{main}. 
As in that case, the remaining task is to find the degree of $ H_n^d $, that is, the dimension of the 
corresponding polytope. Again we have only to find the dimension of the affine subspace of $\R^{n^d}$ 
in which this polytope lives. We do so by counting linearly independent constraint equations. 

To this end, let us write the coordinates of a point in the polytope corresponding to 
$H_n^d$ as $ c(a_1,\dots,a_d) $ with $ 1 \leq a_j \leq n $. These are real 
numbers that satisfy the constraints 
  \[ c(a_1,\dots,a_d) \geq 0 \ , \qquad \sum_{j=1}^n c(a_1,\dots,a_{k-1},j,a_{k+1},\dots,a_d) = 1 \ ( 1 \leq k \leq d ) \ . \]  
Once the $(n-1)^d$ coordinates $ c(a_1,\dots,a_d) $ with $ 1 \leq a_j \leq n-1 $ 
are chosen, the remaining coordinates are clearly determined by the foregoing conditions. Therefore, 
there are at least $n^d - (n-1)^d$ linearly independent constraint equations. 

On the other hand, every constraint involves at least one of the $n^d - (n-1)^d$ coordinates $ c(a_1,\dots,a_d) $ 
for which at least one $a_j$ equals $n$. Consider, say, the coordinate 
$ c(a_1,\dots,a_k,n,\dots,n) $ in which $0 \leq k < d$ and $ a_1,\dots,a_k < n $. We compute: 
\begin{eqnarray*} c(a_1,\dots,a_k,n,\dots,n) &=& 1 - \sum_{ j_d = 1 }^{ n-1 } c(a_1,\dots,a_k,n,\dots,n,j_d) \\ 
                                             &=& 1 - \sum_{ j_d = 1 }^{ n-1 } \left( 1 - \sum_{ j_{d-1} = 1 }^{ n-1 } c(a_1,\dots,a_k,n,\dots,n,j_{d-1},j_d) \right) \\ 
                                             &=& 1 - (n-1) + \sum_{ 1 \leq j_{d-1},j_{d} \leq n-1 } c(a_1,\dots,a_k,n,\dots,n,j_{d-1},j_d) \\ 
                                             &=& \dots \ = 1 - (n-1) + (n-1)^2 - \dots + (-1)^{d-k-1} (n-1)^{d-k-1} \\ 
					     &\mbox{}& \qquad \qquad + (-1)^{ d-k } \sum_{ 1 \leq j_{k+1},\dots,j_{d} \leq n-1 } c(a_1,\dots,a_k,j_{k+1},\dots,j_{d}) \\ 
                                             &=& \frac{ 1 - (1-n)^{d-k} } n + (-1)^{ d-k } \sum_{ 1 \leq j_{k+1},\dots,j_{d} \leq n-1 } c(a_1,\dots,a_k,j_{k+1},\dots,j_{d}) \ . \end{eqnarray*} 
We could have changed the order of summation, that is, used the $d-k$ 
constraints that are in force here in a different order. However, we would 
always end up with the same expression. This implies 
that \emph{all} constraints involving $ c(a_1,\dots,a_k,n,\dots,n) $ with $0 \leq k < d$ and $ a_1,\dots,a_k < n $, 
but not involving any coordinate with $a_j = 1$ for some $j$ with $j \leq k$, are 
equivalent to the \emph{one} constraint 
  \[ c(a_1,\dots,a_k,n,\dots,n) = \frac{ 1 - (1-n)^{d-k} } n + (-1)^{ d-k } 
\sum_{ 1 \leq j_{k+1},\dots,j_{d} \leq n-1 } c(a_1,\dots,a_k,j_{k+1},\dots,j_{d}) \ . \] 
Therefore there are at most $n^d - (n-1)^d$ linearly independent constraint equations. 

Accordingly, the dimension of the polytope and the degree of $ H_n^d $ are $(n-1)^d$. 


\section{CLOSING REMARKS.} 
One big open problem is to determine the periods of our counting functions. The 
evidence gained from our data seems to suggest that the periods increase in some fashion 
with $n$. We believe the following is true: 
\begin{conjecture}. The functions $M_n$, $S_n$, and $P_n$ are not polynomials when $n \geq 5$. \end{conjecture} 
Results in \cite{ahmeddeloerahemmecke} lend further credence to the validity of this conjecture. 
As for semi-magic hypercubes, B{\'o}na showed in \cite{bonamagiccubes} that $H_3^3$ really is a quasi-polynomial, 
not a polynomial. We challenge the reader to prove: 
\begin{conjecture}. The function $H_n^d$ is not a polynomial when $n, d \geq 3$. \end{conjecture} 

\noindent 
{\bf ACKNOWLEDGEMENTS}. We are grateful to Jes\'us DeLoera, Joseph Gallian, Richard Stanley, Thomas 
Zaslavsky, and the referees for helpful comments and corrections on earlier versions of this paper. 


\bibliographystyle{plain}

\begin{thebibliography}{10}

\bibitem{ahmeddeloerahemmecke}
M.~Ahmed, J.~DeLoera, and R.~Hemmecke, Polyhedral cones of magic cubes
  and squares, preprint ({\tt arXiv:math.CO/0201108}) (2002); to appear in: 
  J.~Pach, S.~Basu, M.~Sharir, eds., \emph{Discrete and Computational Geometry---The Goodman-Pollack Festschrift}. 

\bibitem{ananddumirgupta}
H.~Anand, V.~C.~Dumir, and H.~Gupta, A combinatorial
  distribution problem, \emph{Duke Math. J.} \textbf{33} (1966) 757--769. 

\bibitem{andrewsws}
W.~S. Andrews, \emph{Magic Squares and Cubes}, 2nd ed., Dover, New York, 1960. 

\bibitem{beckpixton}
M.~Beck and D.~Pixton, The {E}hrhart polynomial of the
  {B}irkhoff polytope, preprint ({\tt arXiv:math.CO/0202267}) (2002); to appear in \emph{Discrete Comp.~Geom.} 

\bibitem{bonamagiccubes}
M.~B{\'o}na, Sur l'\'enum\'eration des cubes magiques, \emph{C. R.
  Acad. Sci. Paris S\'er. I Math.} \textbf{316} (1993) 633--636. 

\bibitem{ehrhart2}
E.~Ehrhart, Sur un probl\`eme de g\'eom\'etrie diophantienne
  lin\'eaire. {I}{I}. {S}yst\`emes diophantiens lin\'eaires, \emph{J. Reine Angew.
  Math.} \textbf{227} (1967) 25--49. 

\bibitem{ehrhartmagic}
\bysame, Sur les carr\'es magiques, \emph{C. R. Acad. Sci. Paris S\'er. A-B} 
  \textbf{277} (1973) A651--A654. 

\bibitem{macdonald}
I.~G. Macdonald, Polynomials associated with finite cell-complexes, \emph{J.
  London Math. Soc. (2)} \textbf{4} (1971) 181--192. 

\bibitem{macmahon}
P.~A.~MacMahon, \emph{Combinatory Analysis}, Chelsea, New
  York, 1960. 

\bibitem{mcmullenreciprocity}
P.~McMullen, Lattice invariant valuations on rational polytopes, \emph{Arch.
  Math. (Basel)} \textbf{31} (1978/79) 509--516. 

\bibitem{mudgal}
A.~Mudgal, \emph{Counting Magic Squares}, undergraduate thesis, IIT Bombay,
  2002.

\bibitem{paslesfranklin}
P.~C. Pasles, The lost squares of {D}r.\ {F}ranklin: {B}en {F}ranklin's
  missing squares and the secret of the magic circle, \emph{Amer. Math. Monthly} 
  \textbf{108} (2001) 489--511. 

\bibitem{pickoverzen}
C.~A. Pickover, \emph{The Zen of Magic Squares, Circles, and Stars},
  Princeton University Press, Princeton, 2002. 

\bibitem{spencermagic}
J.~Spencer, Counting magic squares, \emph{Amer. Math. Monthly} \textbf{87}
  (1980) 397--399. 

\bibitem{stanleymagic}
R.~P.~Stanley, Linear homogeneous {D}iophantine equations and magic
  labelings of graphs, \emph{Duke Math. J.} \textbf{40} (1973) 607--632. 

\bibitem{ziegler}
G.~M.~Ziegler, \emph{Lectures on Polytopes}, Springer-Verlag, New York,
  1995. 

\end{thebibliography}
\def\cprime{$'$}
\providecommand{\bysame}{\leavevmode\hbox to3em{\hrulefill}\thinspace}
\providecommand{\MR}{\relax\ifhmode\unskip\space\fi MR }
\providecommand{\MRhref}[2]{%
  \href{http://www.ams.org/mathscinet-getitem?mr=#1}{#2}
}
\providecommand{\href}[2]{#2}

\newpage
\setlength{\parindent}{0pt}
\setlength{\parskip}{0.4cm}
\small 
{\bf MATTHIAS BECK} spent his undergraduate years in W\"urzburg, Germany, where he
also enjoyed a short-lived career as a street musician. He got his Ph.D.~with
Sinai Robins at Temple University in 2000 and is currently a postdoc at Binghamton University (SUNY), 
where he enjoys teaching and working with such wonderful students
as Moshe, Jessica, and Paul. His research is in geometric combinatorics and
analytic number theory. Particular interests include counting integer points in all
kinds of polytopes and the application of these enumeration functions to various
combinatorial and number-theoric topics and problems. When he is not counting,
he enjoys biking, travelling, and hanging out with his wife Tendai. \\ 
{\it Department of Mathematical Sciences, Binghamton University (SUNY), Binghamton, NY 13902-6000\\ 
matthias@math.binghamton.edu} 

{\bf MOSHE COHEN} will receive his B.S.~in mathematical sciences at
Binghamton University (SUNY) in May 2004, after which he plans to devote more
time to research in graduate school or an industrial setting.  His
current interests include combinatorial geometry and graph theory.  As
director of an on-campus sound, stage, and lighting business, he can
often be found enjoying concerts from backstage. \\ 
{\it 177 White Plains Rd.~Apt.~8X, Tarrytown, NY 10591 \\ 
bj91859@binghamton.edu} 

{\bf JESSICA CUOMO} is currently pursuing a B.S. in mathematical sciences at Binghamton
University (SUNY) and serves as the president and webmaster of the university's
MAA student chapter. She has more recently studied the field of
combinatorial jump systems at the summer REU at Trinity University and
hopes to continue research in this or similar areas of mathematics.
Her research interests lie mainly in combinatorics, and in the geometry
of discrete systems. \\ 
{\it Dickinson Community \#07120, Binghamton University (SUNY), Binghamton NY 13902 \\ 
jessica@math.binghamton.edu} 

{\bf PAUL GRIBELYUK}, 20, was originally born in Moscow, Russia, and has lived 
in Germany, Arizona, New Jersey, Texas, and New York.  He currently 
attends the University of Texas at Austin, where he is completing majors 
in mathematics and physics.  His current research is in the mathematics of 
finance.  His future plans involve attending graduate school in 
mathematics and working on a topic related to the mathematics of finance. \\ 
{\it 4404 E.~Oltorf St.~\#1202, Austin, TX 78741 \\ 
pavel@math.utexas.edu} 

\end{document}